# Graph-Theoretical Based Algorithms for Structural Optimization

Farzad S. Dizaji[1] Mehrdad S. Dizaji[2]

**Abstract**

Five new algorithms were proposed in order to optimize well conditioning of structural matrices. Along with decreasing the size and duration of analyses, minimizing analytical errors is a critical factor in the optimal computer analysis of skeletal structures. Appropriate matrices with a greater number of zeros (sparse), a well structure, and a well condition are advantageous for this objective. As a result, a problem of optimization with various goals will be addressed. This study seeks to minimize analytical errors such as rounding errors in skeletal structural flexibility matrixes via the use of more consistent and appropriate mathematical methods. These errors become more pronounced in particular designs with ill-suited flexibility matrixes; structures with varying stiffness are a frequent example of this. Due to the usage of weak elements, the flexibility matrix has a large number of non-diagonal terms, resulting in analytical errors. In numerical analysis, the ill-condition of a matrix may be resolved by moving or substituting rows; this study examined the definition and execution of these modifications prior to creating the flexibility matrix. Simple topological and algebraic features have been mostly utilized in this study to find fundamental cycle bases with particular characteristics. In conclusion, appropriately conditioned flexibility matrices are obtained, and analytical errors are reduced accordingly.

## 1. INTRODUCTION

The application of the digital system to structural analysis issues requires the solving of a large number of system of equations of the pattern:

$$Ax = b \qquad (1)$$

This is applicable for either displacement or force techniques. Occasionally, minor changes in matrix A have a significant effect on the solution of Eq. (1). Then we may assert that A is unprepared for this solution. The precision with which Eq. (1) is solved may be affected by matrix A's properties. As a result, it is critical to investigate these features and their connections to the origin, propagation, and spread of potential errors. To do this, more effective problem formulation approaches must be developed, as well as strategies for anticipating, identifying, and reducing solution errors. For the displacement method of analysis, Shah [1] studied the effects of ill conditioning of stiffness matrices. Methods for optimizing the conditioning of stiffness matrices were proposed in his work. Rosanoff and Ginsburg [2] are responsible for a mathematical study of matrix error analysis. They demonstrated that numerically volatile equations might emerge in physically stable situations in their study. Thus, the need to routinely measuring matrix conditioning values associated with different formulation patterns is underlined. Grooms and Rowe [3] examined the impact of sub-structuring on the conditioning of stiffness matrices and found that sub-structuring had no discernible effect on the solution accuracy of ill-conditioned algorithms. Robinson and Haggenmacher investigated how to improve the conditioning of equilibrium equations when using the algebraic force method [4]. Henderson [5] reports that research on the combinatorial force technique has been restricted to improving the sparsity of cycle basis incidence matrices. Recently, techniques for choosing certain kinds of static and kinematical bases have been established, resulting in more constrained stiffness and flexibility matrices compare to traditional ones, Kaveh [6]. In structural engineering, one of the primary causes of ill-conditioning is the usage of components having widely variable

---
[1] Corresponding author, Researcher and Lecturer Professor, Department of Engineering & Society (Applied mathematics and Science Department), University of Virginia, Charlottesville, VA 22901 (Email: ffs5da@virginia.edu), (Phone Number: +1 321 310 5198)

[2] Ph.D. Postdoctoral Research Associate at the Department of Civil Engineering and Environment, University of Virginia, Charlottesville, VA 22901 (Email: ms4qg@virginia.edu), (Phone Number: +1 434 987 9780)

stiffness inside a structure (or flexibility). The use of conventional statical or kinematical basis results in ill-conditioned structural matrices, notwithstanding their optimality. The purpose of this article is to create novel methods for generating special cycles and associated statical and bases that offer the optimal conditioning for flexibility matrices.

## 2. Matrix Force Method

Engineers prefer the force method of structural analysis, in which member forces are utilized as uncertainties since the characteristics of structural members are often determined by element forces instead of just joint deformations. This technique was widely utilized until 1960. Following this, the introduction of the digital computer and the computational feasibility of the displacement technique drew the majority of researchers. Therefore, the force technique and many of the benefits it provides in nonlinear analysis and optimization have been overlooked. The structural analysis force method is divided into four distinct frameworks: 1. Topological force methods, 2. Algebraic force methods, 3. Mixed algebraic-combinatorial force methods, and 4. Integrated force method.

Henderson [5] and Maunder [7] established topological techniques for rigid-jointed skeletal structures by manually selecting the cycle bases of their graph models. Kaveh is credited with developing methods appropriate for computer programming [8-10]. These topological techniques may be used for various kinds of skeletal structures, including rigid-jointed frames, pin-jointed planar trusses, and ball-jointed space trusses [11]. Denke [12], Robinson [13], Topçu [14], Kaneko et al. [15], Soyer and Topçu [16], and Gilbert et al. [17], Coleman and Pothen [18, 19], and Pothen [20] developed algebraic techniques. Patnaik pioneered the integrated force approach [21, 22], in which the equilibrium equations and compatibility criteria are fulfilled concurrently in respect to the force variables. Kaveh et al. [23] used optimization techniques and force method for simultaneous analysis and design of trusses. The graph theoretical force method applied in developing an efficient method for formation of a null bases of finite element models including of tetrahedron, triangle, and rectangular plate elements. Kaveh et al. [24, 25].

Imagine a frame structure $S$ that has $M(S)$ members and $N(S)$ nodes and is statically indeterminate $(S)$ times. Select unknown forces that are $(S)$-independent as redundant. These unknown forces may be selected from the structure's external responses and/or internal forces. Utilize these redundant by

$$q = \{q_1, q_2, ...., q_{\gamma(s)}\}^t \quad (2)$$

To achieve a statically determinate structure, redundant restrictions should be eliminated. This structure is referred to as $S$'s primary (primary or released) structure. This basic structure is considered to be rigid. Consider the external joint loads in the following manner:

$$p = \{p_1, p_2, ...., p_n\}^t \quad (3)$$

Where n is the number of applied nodal force components. For a typical linear analysis using the force method, the stress-resultant distribution owing to a given load $p$ may be expressed as

$$r = B_0 p + B_1 q \quad (4)$$

where $B_0$ and $B_1$ are rectangular matrices, respectively, with $m$ rows and $n$ and $(S)$ columns. The $m$ denotes the number of member forces with independent components. $B_0 p$ is a specific solution that is in equilibrium with the applied stresses, while $B_1 q$ is a complementary solution generated from a maximum collection of self-equilibrating independent stress systems, referred to as a statical basis [26]. While unique and complementary solutions are often derived from the same fundamental framework, this is unnecessary. A basic structure does not have to be chosen arbitrarily. For a redundant basic structure, one may get the required data either by evaluating it first for loads $p$ and bi-actions $q_i = (i = 1, 2, ..., \delta(S))$ or by using pre-existing information [26]. One may write by using the load-displacement relationship for each component and gathering them in the diagonal of the unassembled flexibility matrix Fm, one can write:

$$u = F_m r = F_m B_0 p + F_m B_1 q \quad (5)$$

where $u$ denotes the distortions in the member caused by the internal forces r. The displacements corresponding to the $p$ and $q$ vectors are determined using the contagredient principle.

$$\begin{bmatrix} v_0 \\ v_c \end{bmatrix} = \begin{bmatrix} B_0^t \\ B_1^t \end{bmatrix} [F_m] [B_0 \quad B_1] \begin{bmatrix} p \\ q \end{bmatrix} \tag{6}$$

where $v_0$ represents the displacements associated with the force components of $p$, and $v_c$ represents the relative displacement of the basic structure's released position (cuts). By performing the multiplication, Equation (6) results in

$$\begin{bmatrix} v_0 \\ v_c \end{bmatrix} = \begin{bmatrix} B_0^t F_m B_0 & B_0^t F_m B_1 \\ B_1^t F_m B_0 & B_1^t F m B_1 \end{bmatrix} \begin{bmatrix} p \\ q \end{bmatrix} \tag{7}$$

Imposing the compatibility conditions as $v_c = 0$, the redundant forces are obtained from Equation (7) as

$$q = -(B_1^t F_m B_1)^{-1} (B_1^t F_m B_0) p \tag{8}$$

Substituting in Equation (3), the stress resultant in a structure can be obtained as

$$r = \left[ B_0 - B_1 (B_1^t F_m B_1)^{-1} B_1^t F_m B_0 \right] p \tag{9}$$

In which $G = B_1^t F_m B_1$ is referred to as the structure's flexibility matrix. For a force technique to be efficient, the matrix G should be (a) sparse, (b) well-conditioned, and (c) well-structured, i.e., tightly banded. To ensure that G has the characteristics (a) and (b), the structure of $B_1$ must be precisely constructed since the pattern of $F_m$ remains constant for any discretization; i.e., an appropriate statical basis should be chosen. This issue is addressed in a variety of ways using a variety of techniques. The next sections discuss graph-theoretical techniques for constructing suitable statical bases for various kinds of skeletal systems. The aforementioned characteristic (c) is completely combinatorial in nature. $B_1$ is the pattern comparable to $C^t$ when it has a statical basis in partitioned form. Similar to $CIC^t$ or $CC^t$, $B_1 F_m B_1$ is a pattern that is similar to $CIC^t$ or $CC^t$. This association converts certain structural problems related to defining $G = B_1 F_m B_1$ into combinatorial problems related to dealing with $CC^t$. For instance, if a sparse matrix G is needed, the sparsity of $CC^t$ may be increased. Similarly, with a banded G, one may arrange the corresponding cycles rather than the components of a statical basis (self-equilibrating stress systems). This transformation offers a number of benefits, not the least of which is that the dimension of $CC^t$ is often less than that of G. For example, the dimension of $CC^t$ is six times that of G in a space frame and three times that of G in a planar frame. As a result, when combinatorial characteristics are utilized, the optimization process becomes much easier. Second, since the entries of C and $CC^t$ are elements of $Z_2$, they are simpler to work on than the entries of $B_1$ and G, which are real integers. Thirdly, developments in combinatorial mathematics and graph theory make structural issues immediately relevant. Finally, a connection is established between algebraic and graph-theoretic techniques.

## 3. MINIMAL AND OPTIMAL FUNDAMENTAL GENARALIZED CYCLE BASES

A matrix is said to be sparse if a large proportion of its elements are zero. A fundamental cycle basis $C = \{C_1, C_2, ..., C_{\chi(s)}\}$ is said to be minimal if it corresponds to the smallest value of the following:

$$L(C) = \sum_{i=1}^{b_1(S)} L(C_i) \tag{10}$$

Clearly $\chi(C) = L(C)$ and a minimal GCB may be defined as a basis that corresponds to the smallest $\chi(C)$. A Generalized Cycle Bases (GCB) with a near-minimal $L(C)$ is referred to as a sub-minimal GCB of $S$. The optimum generalized cycle basis of S is a GCB that corresponds to the maximum sparsity of the GCB adjacency matrix. If $\chi(CC^t)$ is not significantly different from its minimal value, the associated basis is said to be suboptimal. The matrix intersection coefficient $\sigma_i(C)$ of row $i$ of the GCB incidence matrix $C$ is equal to the number of rows $j$ such that; $j \in \{i+1, i+2, ... b_1(S)\}$ and $C_i \cap C_j \neq \emptyset$, (i.e., there is at least one $k$ such that the column $k$ of both -cycles $C_i$ and $C_j$ (rows $i$ and $j$) has non-zero entries).

Now it can be shown that:

$$\chi(CC^t) = b_1(S) + 2 \sum_{i=1}^{b_1(S)-1} \sigma_i(C) \tag{11}$$

This connection illustrates the relationship between a GCB incidence matrix C and its GCB adjacency matrix. To minimize $\gamma(D)$, the value of $\sum_{i=1}^{b_1(s)-1} \sigma_i(C)$ should be as little as possible, since $b_1(s)$ is a constant for a given structure S; i.e., γ-cycles with the fewest possible overlaps should be chosen.

# 4. SELECTION OF CYCLE BASES: MATHEMATICAL METHODS

Graphs with cycle bases have a wide variety of applications in a variety of areas of engineering. The quantity of effort required in these applications varies according to the cycle basis used. For certain purposes, a basis with shorter cycles lowers the amount of storage and time required; for others, minimum cycle overlaps are required; for these applications, optimum cycle bases are recommended. The development of minimal and sub-minimal cycle bases is addressed first in this section. Following that, the option of choosing optimum and suboptimal cycle bases is examined. Stepanec [27] proposed minimal cycle bases initially, then Zykov [28] improved on them. Kaveh [9] and Cassell et al. [29]have proposed multiple useful methods for choosing sub-minimal cycle bases. Hubicka and Sysl [30] proposed similar techniques for constructing a graph's minimum cycle basis. Kolasinska [31] discovered an example that contradicted Hubicka and Sysl's method. Kaveh [9] proposed a similar hypothesis for planar graphs; nevertheless, Kaveh and Roosta [32] provided a counter-example. Horton [33] proposed a polynomial-time method for locating graphs' minimum cycle bases. Kaveh and Rahami used algebraic graph theory [34]. The advantages of the algorithms created by various writers are addressed in this part; a technique for selecting minimum cycle bases is provided, as are efficient methods for constructing sub-minimal cycle bases. On a member, the formation of a minimum cycle: A cycle of minimum length $C_i$ on a member $m_j$, referred to as its generator, may be created as follows using the Shortest Route Tree (SRT) algorithm: Begin by forming two SRTs rooted at the $m_j$'s two end nodes $n_s$ and $n_t$, and stop the process when the SRTs intersect at $n_c$ (not via $m_j$). The shortest routes between ns and nc, as well as between $n_t$ and $n_c$, create a minimum cycle $C_i$ on $m_j$. Cycles of specified lengths may also be created using this method. $C_i$, for instance, is a minimum cycle on $m_j$ in Figure 1. SRTs are denoted with bold yellow lines. Take note that the production of SRTs is halted immediately upon the detection of $n_c$ [35].

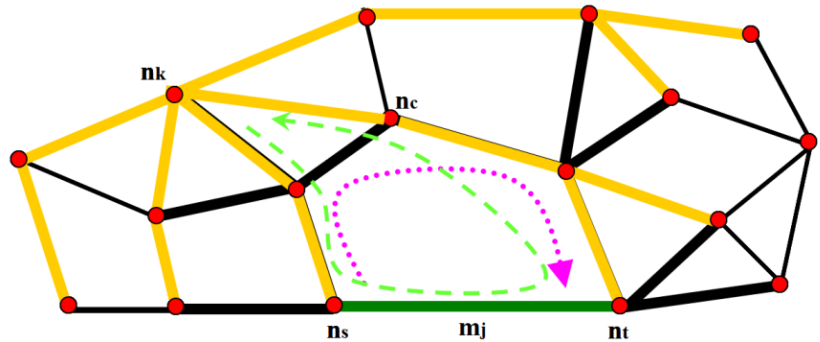

**Figure 1.** A minimal cycle on a graph

Similarly, a minimum cycle on a member $m_j$ that passes through a given node may be produced. An SRT rooted at is constructed, and after it reaches the end nodes of $m_j$, the shortest routes between and $n_s$, and $n_t$ are determined via backtracking. Together with $m_j$, these routes create the necessary cycle. As an example, Figure 1 depicts a minimum cycle on $m_j$, including dashed colored lines. Control of Independence: Each cycle in a graph may be thought of as a column vector in the matrix of its cycle-member incidences. After that, an algebraic technique such as Gaussian elimination may be used to verify a cycle's independence from the previously chosen sub-basis. While this technique is generic and minimizes the order dependence of cycle selection algorithms, its use, like that of many other algebraic methods, needs a significant amount of storage space. The most apparent graph-theoretic method is to construct S's basic cycles using a spanning tree. This technique is effortless; nevertheless, its usage often results in lengthy cycles. The technique may be improved by include each chord in the chosen tree's branch set. Additional length reduction may be accomplished by creating an SRT from a graph's center node and using its chords in increasing order of distance from the center node, Kaveh [8]. A third, similarly graph-theoretical approach is to use admissible cycles. Consider the process of expansion [35].

$$C_1 = C^1 \to C^2 \to C^3 \to .... \to C^b 1^{(S)} = S$$

Where $C^k = \bigcup_{i=1}^{k} c_i$  A cycle C$_{k+1}$ is called an admissible cycle, if for $C^k \rightarrow C^k \cup C_{k+1}$

$$b_1(C^{k+1}) = b_1(C^k \cup C_{k+1}) = b_1(C^k) + 1 \qquad (12)$$

In this paper, methods one and three were used to control the independence of the cycles.

## 5. OPTIMALLY IMPROVED CYCLE BASES

To improve the conditioning of flexibility matrix, it is necessary to choose unique statical bases, or more precisely, cycle bases with specific characteristics. If a cycle basis satisfies the following requirements, it is characterized as an appropriately conditioned cycle basis: (a) It is an optimal cycle basis, i.e., the number of non-zero entries in the corresponding cycle adjacency matrix is minimal, resulting in the flexibility matrix being maximally sparse; (b) The members with the greatest weight of S are involved in the overlaps of the cycles, i.e., the off-diagonal terms of the associated flexibility matrix; and (c) There may be many optimum cycle bases in a weighted graph. The only satisfactory condition (b) is optimally conditioned. If, although, such a cycle basis does not exist, a compromise must be made that satisfies criteria (a) and (b). That is to say, a basis should be chosen that partly meets both requirements. Due to the lack of a method for forming an optimum cycle basis, one should search for a cycle basis with a suboptimal condition.

### 5.1. FORMULATION OF THE PROBLEM

The task of choosing an ideally conditional cycle basis may be expressed mathematically as follows:

$$Min \sum_{i=1}^{b_1(\bar{S})-1} L(C^i \cap C_{i+1}) \qquad\qquad Max \sum_{i=1}^{b_1(\bar{S})-1} W(C^i \cap C_{i+1})$$

Where $\bar{S}$ is a contracted S, and $C^i = \bigcup_{j=1}^{i} C_j$. As can be seen Figure 2, The challenge is a multi-objective optimization one, and the suggested methods are intended to partly satisfy both target functions concurrently.

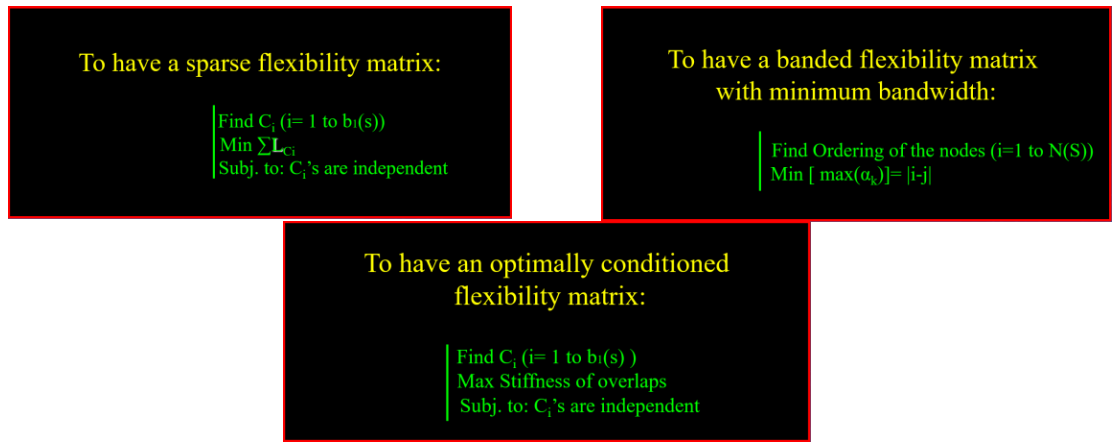

**Figure 2.** Three properties of structural matrices which is a multi-objective optimization problem

## 6. COMPUTATIONAL ERRORS

In numerical modeling, the analysis of errors, either manually or by machine, is essential. Rarely, the input data are precise. Because these data, in most of the cases, are collected based on observations and estimations. Additionally, numerical methods produce different types of errors by themselves. Before studying the topic of errors, the importance of errors was discussed; the goal is to find the roots of the quadratic Equation $x^2 + 0.4002x + 0.00008 = 0$ with four decimals. It turned out that one of the roots of the Equation is -0.0003. The real root with eight decimals is -0.0002. Therefore, calculus with four decimals produces 50% error. The reason for this error is using the calculus with four-digit decimals. Consequently, we should not count on eight decimals and assume that eight-digit decimals are the solution here.

**6.1. Ill-conditioning and instability.**

Some of the calculations are very sensitive to round-off errors; however, some do not. In some cases, this sensitivity could be addressed with the change in the analysis method, but in most cases, the defined problem is intrinsically sensitive to the round-off errors. Therefore, sensitivity to the method and intrinsic sensitivity must be distinguished and treated separately. Generally, the word of sensitivity can be described as a continuous dependency of a problem on the input data and its method. As in a problem, if the rate of an error increases continuously, the method will be called unstable numerically. As an example, we are Investigating the ill-conditioning of $(x - 2)^4 = 0$; this Equation has four same roots, which is x = 2. Now, subtract the same Equation from a minimal number; for instance, $10^{-8}$. Then study the behavior of the new Equation.

$$(x-2)^4 - 10^{-8} = x^4 - 4x^3 + 8x^2 - 16x0 + 15.99999999 = 0$$

The roots of the above Equation should satisfy $(x - 2)^4 = \pm 10^{-4}$ which is $(x - 2) = \pm\sqrt{\pm 10^{-4}}$. Therefore the roots of the Equation are

$$x_1 = 2.01 \quad , \quad x_2 = 1.99 \quad , \quad x_3 = 2+0.01i, \quad x_4 = 2-0.01i$$

As can be seen, the roots of the considered Equation are very sensitive to the input data (e.g., a coefficient). This type of behavior, independent from the solution method and belongs to intrinsic properties of the problem or Equation, leads to ill-conditioning of the whole system. As observed, the roots of the Equation faced significant changes just with changes in the constant coefficient. As the second example for ill-condition systems we are solving the following linear system; $Ax = b$; using Gauss elimination method and considering four digit decimals as round off errors (after each operation just keep four decimals), the solution is

$$[A] = \begin{bmatrix} -0.002 & 4.0 & 4.0 \\ -2 & 2.906 & -5.38 \\ 3.0 & -4.301 & -3.112 \end{bmatrix}, \quad [b] = [7.998 \quad -4.481 \quad -4.143], \quad x = [-1496 \quad 2.0 \quad 0.0]$$

But the precise solution is

$$x = [-1.0 \quad 1.0 \quad 1.0]$$

Matrix A is ill-condition and in order to find the precise solution we need to make some changes to the matrix. In doing so, the rows of the matrix were changed as follow;

$$[A] = \begin{bmatrix} 3.0 & -4.301 & -3.112 \\ -0.0002 & 4.0 & 4.0 \\ -2.0 & 2.406 & -5.386 \end{bmatrix}, \quad [b] = [-4.413 \quad 7.998 \quad -4.481], \quad x = [-1.0 \quad 1.0 \quad 1.0]$$

With changing the rows of the equations (matrix), a big difference observed in the solution. This kind of systems called ill-condition systems.

Generally, a matrix is ill-condition if at least one of the following conditions satisfies:
1. A small change in coefficients (arrays of the matrix) causes big changes in the solution of the linear equations
2. Diagonal arrays are much smaller than the rest of the arrays.
3. The value of the $|A|.|A^{-1}|$ never becomes one.
4. The value of the $[A^{-1}]^{-1}$ is not precisely the same as A
5. The value of the $[AA^{-1}]$ is not exactly equal to I (identity matrix)

**6.2. CONDITION NUMBERS**

Numerous numbers are established and used in practice to quantify the conditioning of a matrix. Three frequently used condition numbers are described here; they are straightforward and simple to use.

### 6.2.1. THE RATIO OF EXTREME EIGENVALUES

Eigenvalues and eigenvectors are closely linked to matrix conditioning. The ratio of a matrix $|\lambda_{max}|/|\lambda_{min}|$ extreme eigenvalues may be used to determine its condition number. It is straightforward to demonstrate that the logarithm of this condition number to base ten is approximately proportional to the maximum number of significant figures lost during inversion or solution of simultaneous equations. Thus the number of good digits, g, in the solution is given by:

$$g = p - \log(|\lambda_{max}|/|\lambda_{min}|) = p - PL$$

In this correlation, $PL = \log(|\lambda_{max}|/|\lambda_{min}|)$ and p is a variable that changes across machines. For example, the IBM/360 does single-precision computations with about eight digits and double-precision calculations with around sixteen digits. It must be emphasized that the estimate above is conservative, and history suggests that PL should be on the conservative side by one digit., Kaveh [35].

### 6.2.2. DETERMINANT OF A ROW-NORMALIZED MATRIX

A straightforward and practicable way to determine the conditioning of a set of equations is to compute the determinant of the row-normalized matrix of the set's coefficients. This implies that each row of matrix A, let us call Ai, is split by the following:

$$\left[ g_{i1}^2 + g_{i2}^2 + \ldots + g_{ii}^2 \right]^{1/2}$$

The size of the row-normalized A's determinant, indicated by PN, is a useful indicator of A's conditioning. This determinant's size falls inside the range,

$$0 < PN \leq 1$$

Because A is always positive definite. In the case of orthogonal or diagonal matrices, the matrix with perfect conditioning has PN = 1.

### 6.2.3. WEIGHTED GRAPH AND AN ADMISSIBLE MEMBER

The relative stiffness (or flexibility) of structure members may be regarded as positive integers associated to the structure's graph model, resulting in a graph which is weighted. Assume *S* represent the model of a frame structure, and $k_{mi}$ is the stiffness matrix of an element $m_i$ in the structure's global coordinate system. Using the diagonal elements $k_{ii}$ of $k_{mi}$, a weight may be determined for $m_i$., as

$$W(m_i) = \sum k_{ii} = 2(\alpha_1 + \alpha_4^z + \alpha_3^z)$$

$$\alpha_1 = \frac{EA}{L} \quad , \quad \alpha_4^z = \frac{12EI}{L^3} \quad , \quad \alpha_3^z = \frac{4EI}{L}$$

Additionally, an alternative weight based on the square roots of $k_{mi}$'s diagonal entries may be used:

$$W(m_i) = \sum k_{ii}^{1/2} = 2\left[ (\alpha_1)^{1/2} + (\alpha_4^z)^{1/2} + (\alpha_3^z)^{1/2} \right]$$

Other weights functions may be constructed as appropriate to reflect the relative stiffness of the members of *S*. Definition: Assume the weight of members $m_1, m_2... m_{M(S)}$ be defined by $W(m_1), W(m_2)…W(m_{M(S)})$ respectively. A member mi is called F-admissible, if

$$W(m_i) = \frac{1}{\alpha} \sum_{j=1}^{M(s)} W(m_j)/M(S)$$

Where α is an integer number which can be taken as 2, 3, ... We utilized α = 2; however, a thorough research using different quantities of α is necessary. Inadmissible or S-admissible members are those who are not F-admissible. Figure 3 illustrates two distinct structural models and their graph-based counterparts.

## 7. BASIC DEFINITIONS AND CONCEPTS f GRAPH THEORY

A graph S is made up of a set N(S) of nodes (vertices or points) and a set M(S) of members (edges or arcs), as well as an incidence relation for each member of a pair of distinct nodes, referred to as its ends. A skeletal structure's connectivity characteristics may conveniently be translated to those of a network S; the structure's joints and members correspond to the graph's nodes and edges, respectively.

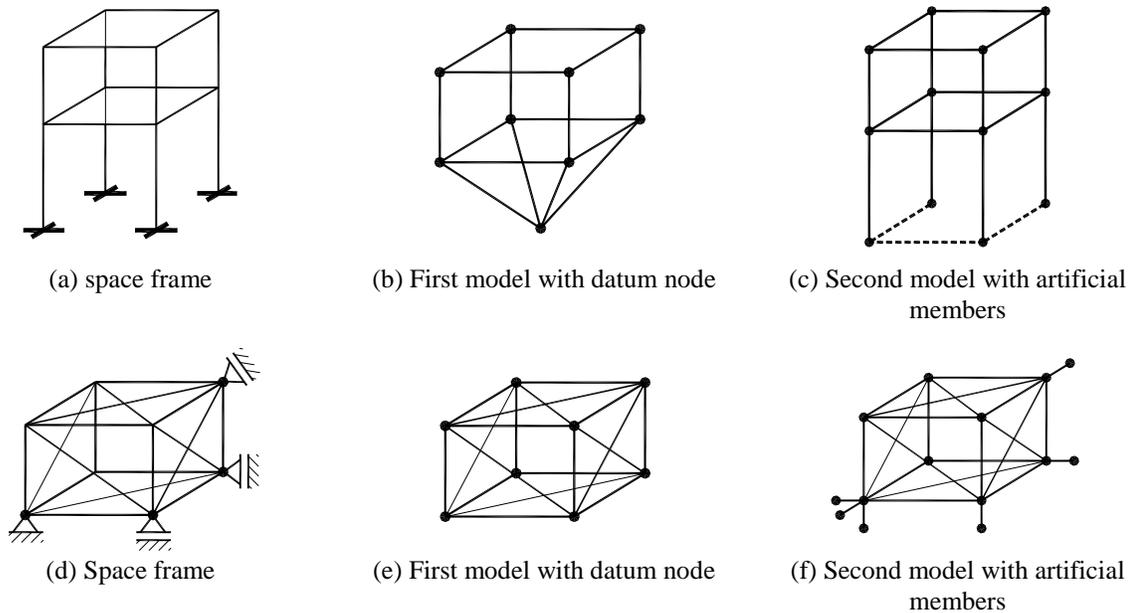

(a) space frame     (b) First model with datum node     (c) Second model with artificial members

(d) Space frame     (e) First model with datum node     (f) Second model with artificial members

**Figure 3.** Various structural model and their associated graph model

If a network (graph) can be plotted on a plane and its members cross only at nodes, it is said to be planar. A path across *S* is a finite sequence $P_k = (n_0, m_1, n_1, ..., m_p, n_p)$ whose terms alternate between distinct nodes $n_i$ and members $m_i$ of *S*, where $n_{i-1}$ and $n_i$ are the two ends of $m_i$. The length of a path $P_i$, represented by $L(P_i)$, is used to determine its member. If there is an alternative route $P_j$ between the two end nodes $n_0$ and $n_p$, that satisfies $L(P_i) \le L(P_j)$ then Pi is considered the shortest path between these nodes. The distance between two nodes $n_1$ and $n_2$ in a network S is defined as the length of the shortest route between them, given that the path's elements are included inside S, and is denoted by $d_S(n_1, n_2)$ or $d_S(n_2, n_3)$. In S, two nodes $n_i$ and $n_j$ are said to be linked if a route exists between them. If all pairs of nodes in a network S are linked, the graph is said to be connected. Components of S are maximally connected subgraphs, i.e., they are not subgraphs of any other connected subgraph of S. A graph is said to be two-connected if it retains its connectivity after the removal of one of its members. A cycle is a route that contains $(n_0, m_1, n_1, ..., m_p, n_p)$ for which $n_0 = n_p$ and $p \ge 1$; in other words, a cycle is a closed path. T of S is a linked subgraph that does not include any cycles. A tree is said to be a spanning tree of S if it includes all of S's nodes. A shortest route tree (SRT) with its root at a given node $n_0$ of S is a spanning tree in which the distance between each node $n_i$ of T and $n_0$ is as little as possible. The following basic method may be used to create an SRT rooted at a given node $n_0$: Label the chosen root $n_0$ with a value of '0' and the surrounding nodes with a value of '1'. As tree members, set the members' incident to '0'. Repeat the process of labeling with '2' the unnumbered endpoints of all the members that intersect with nodes labeled with '1', recording the tree members once more. This procedure is complete after each node in S has been tagged and all tree members have been documented. Let us define a cycle set of members of a graph as a collection of members that form a cycle or multiple cycles with no common member but maybe shared nodes. A cycle set vector is a vector that represents a cycle set. It can be shown that the sum of a graph's two cycle set vectors is also a cycle set vector. Thus, a graph's cycle set vectors define a vector space over the integer modulo 2 field. A cycle space has a dimension equal to the first Betti number of graph $b_1(S)$.

## 8. FUNDAMENTAL CYCLE BASES

A basic cycle basis, a kind of special cycle basis, maybe readily built to match to a tree T of S. Within a linked S, a chord of T and T includes a cycle referred to as the basic cycle of S. Additionally, the basic cycles formed by adding the chords to T, one at a time, are independent since each cycle has a component that is not present in the others. Additionally, each cycle $C_i$ is dependent on the set of basic cycles produced via the preceding procedure since $C_i$ is the symmetric difference of the cycles defined by the chords of T that fall in q. Thus the cycle rank (cyclomatic number, first Betti number, nullity) of graph S, which is the number of cycles on the basis of the cycle space of S, is given by,

$$rank(S) = b_1(S) = M(S) - N(S) + b_0(S)$$

## 9. THE PROPOSED ALGORITHMS FOR GENERATING SUB-OPTIMAL CYCLE BASES

In this section, five new algorithms were proposed for optimal fundamental cycle bases of a weighted graph. In each selected cycle basis, for a corresponding graph model, there are three or six (S.E.Ss) for 2D or 3D, respectively [36]. For all five algorithms independency of the cycle bases checked using an algebraic method Gaussian elimination to verify a cycle's independence from the previously chosen sub-basis

*Algorithm 1 (this algorithm generates sparser fundamental cycle bases and less well-conditioned)*

- **Step 1.** Create graph model of the structure.

- **Step 2:** From a member of the graph model, create a cycle with minimum length and maximum weight.

- **Step 3:** Repeat step 1 for the rest of the members of the graph model; (i.e., producing a cycle base with minimum length and maximum weight from all members of the graph model. Construct a SRT with maximum weight from two ends of a graph member and return to two ends of the same member with a Path as soon as two SRT cut each other or they meet each other at one node. In doing so, a cycle with minimum length and maximum weight will be created). The shortest route tree (SRT) rooted at a specified node $n_o$ of S is a tree for which the distance between every node $n_j$ of T and $n_o$ is a minimum. For example SRT created from $n_o$ to $n_{15}$ in the Figure 4 which shows the minimum path between two points. The following basic method may be used to create an SRT for a graph: Label the chosen root with the letter "o" and the surrounding nodes with the number "1". Keep track of the members' encounters with "o" as tree members. Repeat the procedure of labeling the unnumbered endpoints of all the members that intersect with nodes labeled as "1" with "2," once more recording the tree members. This procedure is complete after each node in S has been tagged and all tree members have been documented. The label of the last node represents the SRT's length, whereas the maximum number of nodes with the same label specifies the SRT's breadth.

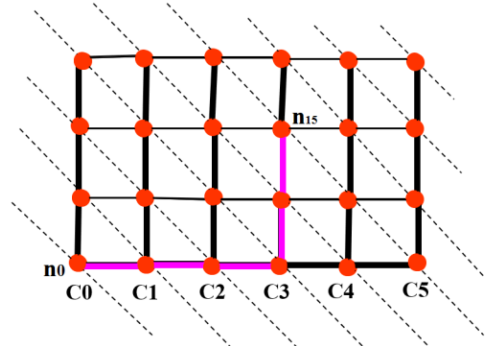

**Figure 4.** SRT created between two points n0 and $n_{15}$

- **Step 4:** After we are done with step 2, a cycle with minimum length and maximum weight was generated from all members of the graph model. The weight of each cycle should be calculated (the weight of a cycle can be obtained by adding all weights of the members of that specific cycle).

- **Step 5:** Putting all created cycles in a weight ascending order. In fact, we are sorting out all cycles according to their weight.

- **Step 6:** The greedy algorithm should be utilized here. We will continue selecting cycle bases until we have $b_1(s)$ independent bases, which will be a set of our fundamental cycle bases (first, using a greedy algorithm, a cycle with maximum weight, from the set of step 4, will be selected, then the next cycle should be selected from the remaining set of step 4 which is independent than first selected cycle base and step 5 should be carried on until we have $b_1(s)$ cycle bases. In doing so, a set of optimal fundamental cycle bases will be generated).

*Algorithm 2* (this algorithm generates improved well-conditioned fundamental cycle bases and less sparser)

- **Step 1:** Create a graph representation of the structure.

- **Step 2:** First, find the graph model of the structure, then generate a cycle with maximum weight and minimum length.

- **Step 3:** Repeat step 1 for the rest of the members of the graph model, producing a cycle base with minimum length and maximum weight from all members of the graph model. (The intention of the previous sentence is that we are compromising between being the minimum length of the cycle and being maximum of the same cycle, which is an objective that cannot be achieved simultaneously. As mentioned before, finding an algorithm that can produce a sparse matrix with excellent well-conditioning is not possible was shown that this is a multi-objective optimization problem and we must compromise between two properties which are sparsity and well conditioning. It means that we should apply sparsity algorithms if the sparsity of the matrix is essential for our problem that makes the matrix sparser and vice versa). Construct a SRT with maximum weight from two ends of a graph member and return to two ends of the same member with a Path as soon as two SRT cut each other or they meet each other at one node. In doing so, a cycle with minimum length and maximum weight will be created.

    Generating the Shortest Root Tree Maximum (SRTM): label the selected note for SRTM as "0" then identify the members who have an incident with the same note. Then calculate the average weight of the members that incident with node "0". Remove the members (incident members with node "0") from the branches of SRTM with a weight less than the average weight. In the process, assign "1" to the end of members that have an incident with node "0". Then identify all members that have an incident with node "1" and sort them as ascending in their weights. Afterward, select the tree members from the ordered set and assign "2" to the opposite end of

members. Next, we obtain the average weights of the members with an incident with node "1" and remove members with less than average weight from the SRTM branch. For the last part, repeat this procedure for node "2" and the rest of the graph model nodes until SRTM is generated.

- **Step 4:** Repeat step 3 for the rest of the members of the graph model. We will generate a cycle with minimum length and maximum weight for all members of the graph. Then the weight of each of the cycles should be calculated (the weight of each cycle can be obtained by summing all weights of the members of that cycle).

- **Step 5:** sort out all cycles as ascending weight and put them in a set.

- **Step 6:** The greedy algorithm should be utilized here. We will continue selecting cycle bases until we have $b_1(s)$ independent bases, which will be a set of our fundamental cycle bases (first, using a greedy algorithm, a cycle with maximum weight, from the set of step 4, will be selected, then the next cycle should be selected from the remaining set of step 4 which is independent than first selected cycle base and step 5 should be carried on until we have $b_1(s)$ cycle bases. In doing so, a set of optimal fundamental cycle bases will be generated).

***Algorithm 3*** *(this algorithm generates sparser fundamental cycle bases and less well-conditioned)*

- **Step 1:** Create a graph representation of the structure.

- **Step 2:** Create a cycle from a member of the graph model that has the shortest possible length and the greatest possible weight.

- **Step 3:** Repeat step 1 for the remaining members of the graph model; (i.e., create a cycle base with a minimum length and maximum weight using all members of the graph model. Construct an SRT with the greatest weight between two ends of a graph member and return to the same member's two ends through a Path as soon as two SRT intersect or meet at a node. Thus, a cycle with the shortest possible cycle and the heaviest possible weight will be produced).

- **Step 4:** After completing step 2, we created a cycle with a minimum length and maximum weight using all members of the graph model. Each cycle's weight should be determined (the weight of a cycle can be obtained by adding all weights of the members of that specific cycle and the length of a cycle can be calculated by adding the length of all members of that cycle).

- **Step 5:** Arranging all generated cycles in increasing length order. Indeed, we are classifying all cycles based on their length.

- **Step 6:** Here, the greedy method should be used. We will continue selecting cycle bases until we have $b_1(s)$ independent bases, which will be a set of our fundamental cycle bases (first, using a greedy algorithm, the cycle with the minimum length should be selected from the set of step 4, then the next cycle should be selected from the remaining set of step 4 that is independent of the first selected cycle base, and so on until we have $b_1(s)$. This generates a collection of optimum basic cycle bases).

***Algorithm 4*** *(this algorithm generates improved well-conditioned fundamental cycle bases and less sparser)*

- **Step 1:** Create a graph representation of the structure.

- **Step 2:** Obtain the lowest possible length and the highest possible weight from a member of the graph model by constructing a cycle from that component. This is almost similar to algorithms but in step 2 we replaced SRT with SRTM.

- **Step 3:** To construct a cycle base with a minimum length and maximum weight using all of the members of the graph model, repeat step 1 for the remaining members of the graph model. Additionally, we need to compromise between finding a cycle with minimum length and maximum weight. Because it is impossible to find a cycle that satisfy both requirement due to the fact that this problem and our algorithms belong to NP categories which means that they don't have one objection function.

- Generating the Shortest Root Tree Maximum (SRTM): label the selected note for SRTM as "0" then identify the members who have an incident with the same note. Then calculate the average weight of the members that incident with node "0". Remove the members (incident members with node "0") from the branches of SRTM with a weight less than the average weight. In the process, assign "1" to the end of members that have an incident with node "0". Then identify all members that have an incident with node "1" and sort them as ascending in their weights. Afterward, select the tree members from the ordered set and assign "2" to the opposite end of members. Next, we obtain the average weights of the members with an incident with node "1" and remove members with less than average weight from the SRTM branch. For the last part, repeat this procedure for node "2" and the rest of the graph model nodes until SRTM is generated.

- **Step 4:** After completing step 2, we used all of the components of the graph model to construct a cycle with a minimum length and a maximum weight. It is necessary to establish the weight of each cycle (the weight of a cycle can be obtained by adding all weights of the members of that specific cycle and the length of a cycle can be calculated by adding the length of all members of that cycle).

- **Step 5:** Putting all of the produced cycles in a sequence of increasing length. In fact, we are categorizing all cycles according to their length.

- **Step 6:** Here, the greedy method should be used. We will continue selecting cycle bases until we have $b_1(s)$ independent bases, which will be a set of our fundamental cycle bases (first, using a greedy algorithm, the cycle with the minimum length should be selected from the set of step 4, then the next cycle should be selected from the remaining set of step 4 that is independent of the first selected cycle base, and so on until we have $b_1(s)$. This generates a collection of optimum basic cycle bases).

*Algorithm 5* (*this algorithm generates improved well-conditioned fundamental cycle bases and less sparser*)

- **Step 1:** Create a graph representation of the structure.

- **Step 2:** Find the weight of all members of the created graph and sort them either in a descending or ascending format, and put them in a set

- **Step 3:** From the set in step 2, select a member with minimum weight and find out whether or not the member is admissible or not. If the selected member is not admissible, put it in another set and name the set *NA*.

- **Step 4:** Select the next member with minimum weight from the set in step 2 and observe its admissibility. If the selected member is not admissible, then out the member in the *NA* list. This step will be continued until we find the first admissible member. This is a list of all inadmissible members. The goal of the prepared list is that we avoid putting them in overlaps of fundamental cycles.

- **Step 5:** We use the member list *NA* again in such a way that the first member, $m_1^g$ with the lowest weight will be selected and generate a cycle with minimum length and maximum weight. (Construct an SRT with the highest weight between the two ends of a graph member and return to the two ends of the same member through a Path whenever two SRT collide or meet at a node. Thus, a cycle with the shortest possible cycle and the highest weight possible will be created).

- **Step 6:** From the *NA* member list the next member with minimum weight will be selected and generate a cycle with minimum length and maximum weight in such a way that does not include members $m_1^g$, which means that none of the members of the next cycle include member $m_1^g$.

- **Step 7:** Repeat step 5 for the rest of the members included in list *NA*. For instance, construct j-cycle with the minimum length and maximum weight, which includes the edges of $S / \bigcup_{i=1}^{j-1} m_i^g$. In fact, the goal is that these members do not contain in the overlap of the cycles. Because their weights are small compare to others members 'weight and this will cause the ill-conditioning of the matrix.

- **Step 8:** As soon as we repeated step 7 for all members in list *NA*, we select a member of the graph model which is not in list *NA* and generate a cycle with minimum length and maximum weight on this selected member.

- **Step 9:** Repeat step 8 for the rest of the members who are not in list *NA* and put them in a set.

- **Step 10:** Sort out all of the cycles in an ascending weight or ascending length and again put them in a set.

- **Step 11:** Here, the greedy method should be used. We will continue selecting cycle bases until we have $b_1(s)$ independent bases, which will be a set of our fundamental cycle bases.

In this algorithm, if in step 9 we sort out cycles in terms of their weight, then use the greedy algorithm, in fact; we prefer cycles with the length of longer, which makes the flexibility matrix less sparse, but in return, the flexibility matrix will be well conditioned and vice versa.

## 10. ILLIRASTRATION EXAMPLES AND RESULTS

In order to investigate the efficiency of the proposed algorithms, a JAVA programming language was written for all five algorithms and a friendly GUI software tool is developed to visualize optimal cycle bases. Moreover, the well conditioning of the five proposed algorithms was compared to each other. The visual form of the programming is shown in Figure 5. Additionally, Figure 6 shows how we are selecting different weights and algorithms in order to find optimal cycle bases. Additionally, in the examples two different set of structural properties were used for the frames. Using structural properties we are able to calculate the weight of each member. Moreover, in the examples two different colors were shown in the frames which represent different properties. For example, black color's properties shown on the frames share same structural properties which is different than the properties depicted with red color. Moreover, fictional yellow color members were attached to the supports of the frames in order to guarantee the rigidity of the graph model [37].

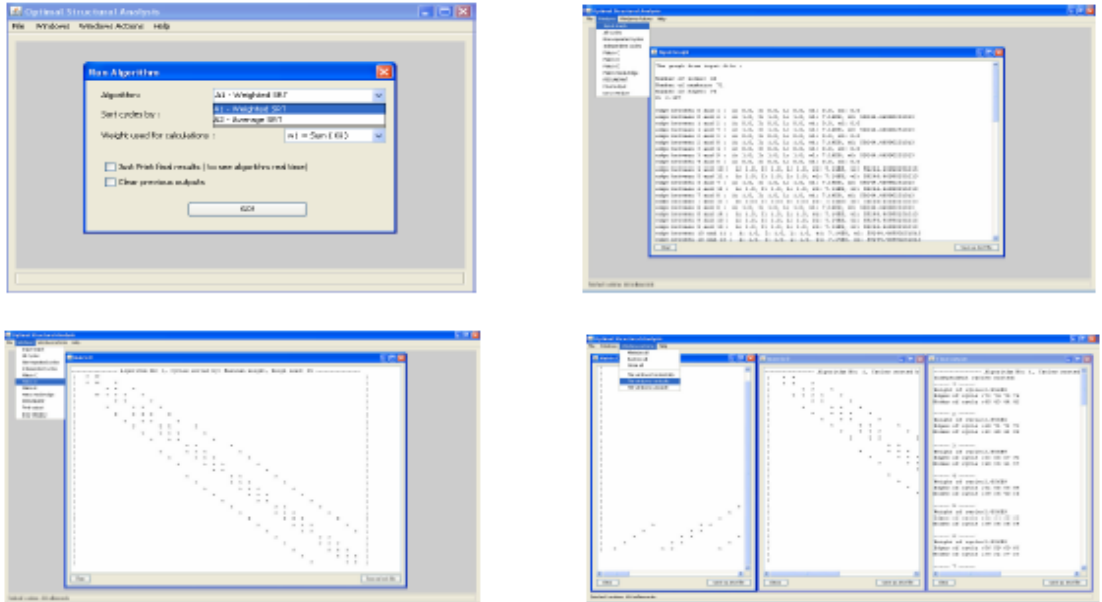

**Figure 5.** The schematic of the written programming in JAVA

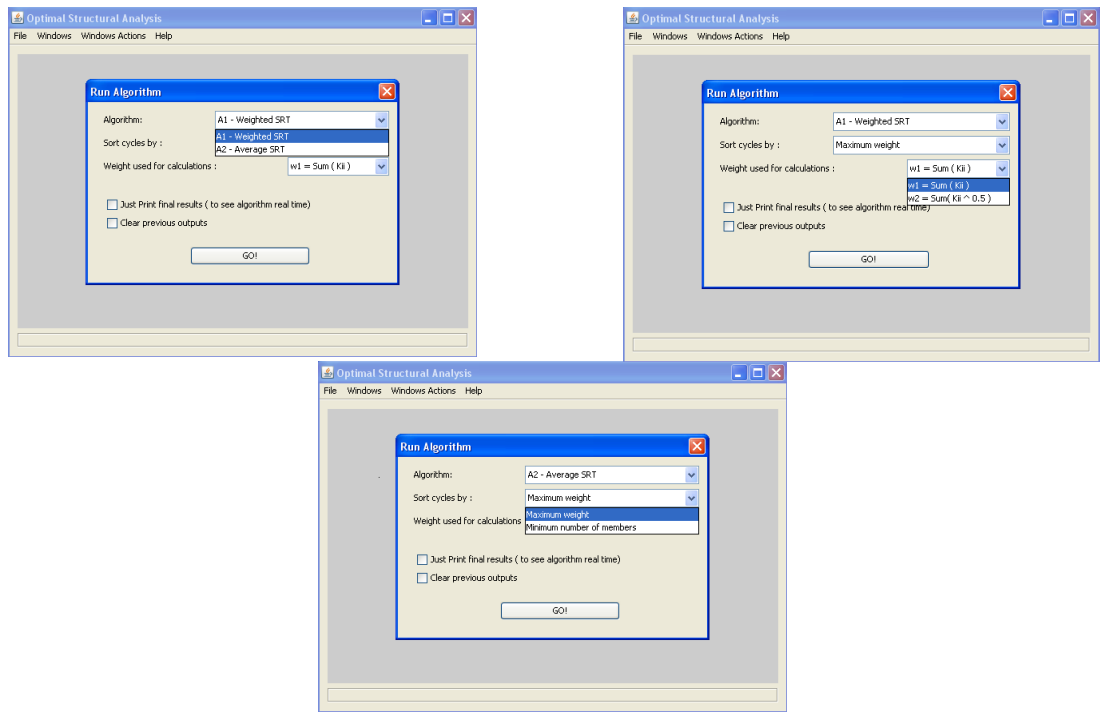

**Figure 6.** Schematic of the selecting the type of algorithm in finding cycle bases

**Example 1** – A 3-story, 4-span frame was selected in Figure 7. The properties of beams and columns were shown in Table 1. Three well conditioning numbers were studied for two members' properties as follows:

**Table 1.** Defined properties for beams and columns of the frame

| Structural Properties | $A_1$ (m$^2$) | $A_2$ (m$^2$) | $I_1$ (m$^4$) | $I_2$ (m$^4$) | E (t/m$^2$) |
|---|---|---|---|---|---|
| | 0.00106 | 0.00970 | 0.00000171 | 0.00019610 | $2.1 \times 10^7$ |

The length of all members (beams and columns) was assigned to be 3m.

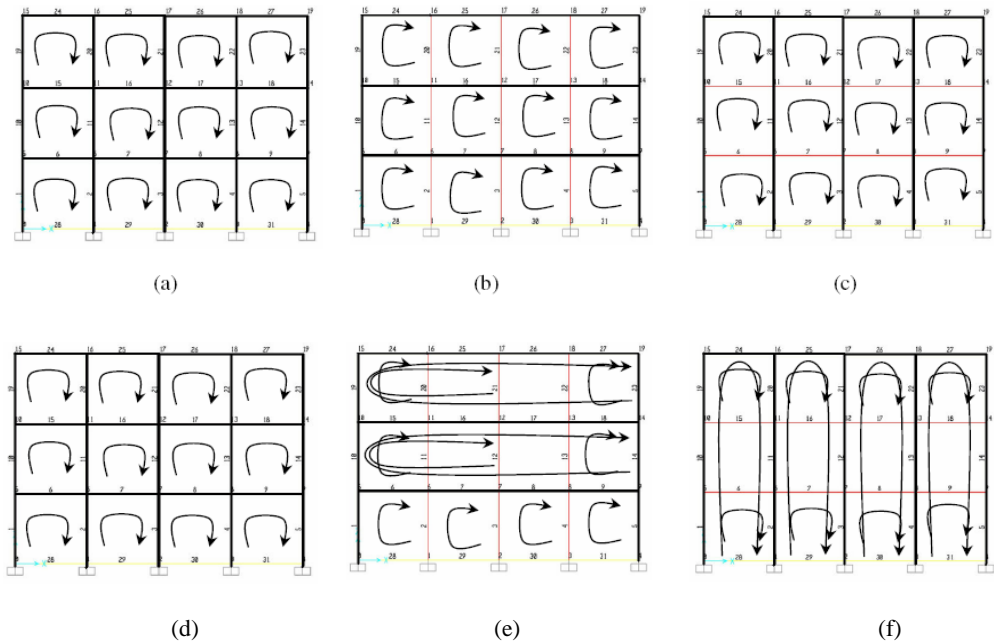

**Figure 7.** 3-story, 4-span frames with different member properties and depicted optimal cycle bases using algorithm 1 (a-c) and algorithm 2 (d-f)

As can be seen in Figure 7, optimal cycle bases were shown on the frames using algorithms 1 and 2. The corresponding flexibility matrices for frames (d-f) were presented in Figure 8. For the frames (a and d), since the same properties were used for beams and columns, both algorithms produced the same length of cycle bases. Also, each of the cycles has the minimum length and that is why the corresponding flexibility matrix is sparser in Figure 8a. On the contrary, since the properties of the beams and columns in the real world are not similar to each other, for the frames (b-c and e-f) various structural properties were used for the beams and columns. Thus, as seen in Figure 8(b-c), corresponding flexibility matrix is less sparse.

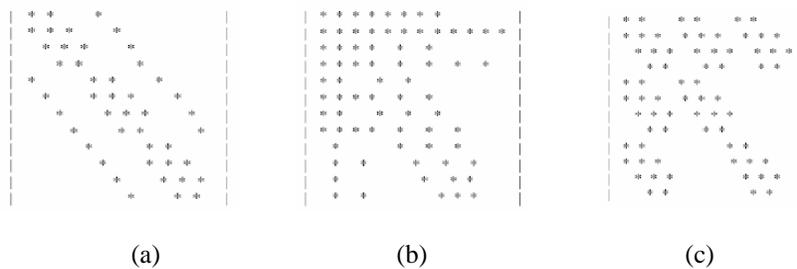

**Figure 8.** Corresponding flexibility matrices for frames (d-f) of Figure 7.

Applying algorithm 1 for the frame with different properties produces sparser matrix compare to algorithm 2. On the other hand, algorithm 2 generates less sparse matrix but with better well conditioning. This fact can clearly be seen, by comparing the conditioning number, in the Table 2 and Table 3.

**Table 2.** Well conditioning numbers using algorithm 1

| TYPE | PL | PN | PDET | X(D) |
|---|---|---|---|---|
| (a) | 3.452154 | 6.381245E-29 | 7.134199E-20 | 46 |
| (b) | 4.852145 | 0 | 0 | 46 |
| (c) | 4.479986 | 0 | 0 | 46 |

**Table 3.** Well conditioning numbers using algorithm 2

| TYPE | PL | PN | PDET | X(D) |
|---|---|---|---|---|
| (d) | 3.452154 | 6.381245E-29 | 7.134199E-20 | 46 |
| (e) | 4.342158 | 2.063098 E-42 | 3.613358E-31 | 74 |
| (f) | 3.985745 | 8.307591E-41 | 3.737124E-30 | 70 |

**Example 2** – A 3-story, 3-span frame was selected in Figure 9. The properties of beams and columns were shown in Table 4. The optimal fundamental cycle bases for the frames generated using algorithms 3 and 4, three well-conditioning numbers were studied for two members' properties as follows:

**Table 4.** Defined properties for beams and columns of the frame

| Structural Properties | $A_1$ (m$^2$) | $A_2$ (m$^2$) | $I_1$ (m$^4$) | $I_2$ (m$^4$) | E (t/m$^2$) |
|---|---|---|---|---|---|
| | 0.00106 | 0.00970 | 0.00000171 | 0.00019610 | $2.1\times 10^7$ |

The length of the all members were assigned to be 3m.

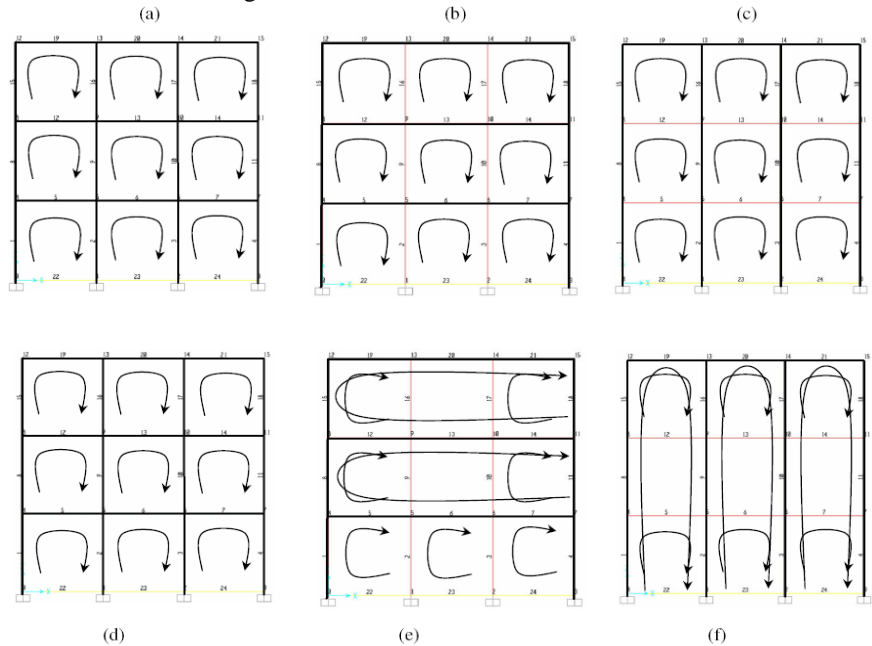

**Figure 9.** 3-story, 3-span frames with different member properties and depicted optimal cycle bases using algorithm 3 (a-c) and algorithm 4 (d-f)

As shown in Figure 9 , optimum cycle bases for the frames were determined using algorithms 3 and 4. Figure 10 depicts the appropriate flexibility matrices for frames (d-f). Due to the fact that the identical characteristics were utilized for beams and columns in both frames (a and d), both algorithms generated the same length of cycle bases. Additionally, each cycle has the shortest possible length, which results in a sparser flexibility matrix. Figure 10a.

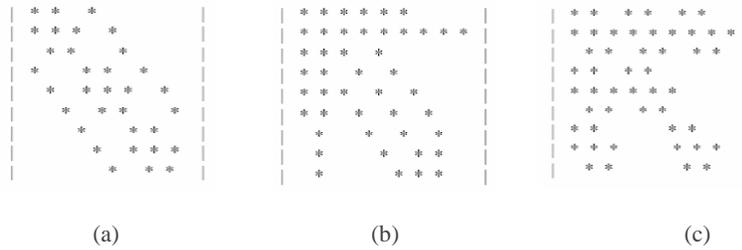

(a)          (b)          (c)

**Figure 10.** Corresponding flexibility matrices for frames (d-f) of Figure 9.

On the contrary, since the characteristics of beams and columns in the actual world are not identical, different structural qualities were utilized for the beams and columns in the frames (b-c and e-f). As shown in Figure 10(b-c), the associated flexibility matrix is thus less sparse. When applied to a frame with differing characteristics, method 3 generates a sparser matrix than algorithm 4. On the other side, method 4 produces a matrix that is less sparse but has improved well conditioning. This fact is readily apparent when the conditioning number in Table 5 is compared to the conditioning number in Table 6.

**Table 5.** Well conditioning numbers using algorithm 3

| TYPE | PL | PN | PDET | X(D) |
|---|---|---|---|---|
| (a) | 3.172154 | 5.732461E-21 | 2.214527E-14 | 33 |
| (b) | 4.632548 | 0 | 0 | 33 |
| (c) | 4.320656 | 6.918279E-40 | 3.524312E-31 | 33 |

**Table 6.** Well conditioning numbers using algorithm 4

| TYPE | PL | PN | PDET | X(D) |
|---|---|---|---|---|
| (d) | 3.172154 | 5.732461E-21 | 2.214527E-14 | 33 |
| (e) | 4.125478 | 3.245248E-34 | 4.572541E-26 | 45 |
| (f) | 3.826745 | 9.192013E-30 | 1.684895E-21 | 49 |

**Example 3** – A 4-story 4-span frame with members with different properties was shown in Figure 11. In order to find sub-optimal fundamental cycle bases algorithms 1 and 2 were applied respectively. Also, the defined properties for beams and columns of the frame were shown in Table 7.

**Table 7.** Defined properties for beams and columns of the frame

| Structural Properties | $A_1$ (m$^2$) | $A_2$ (m$^2$) | $I_1$ (m$^4$) | $I_2$ (m$^4$) | E (t/m$^2$) |
|---|---|---|---|---|---|
| | 0.00106 | 0.00970 | 0.00000171 | 0.00019610 | $2.1 \times 10^7$ |

The length of the all members were assigned to be 3m.

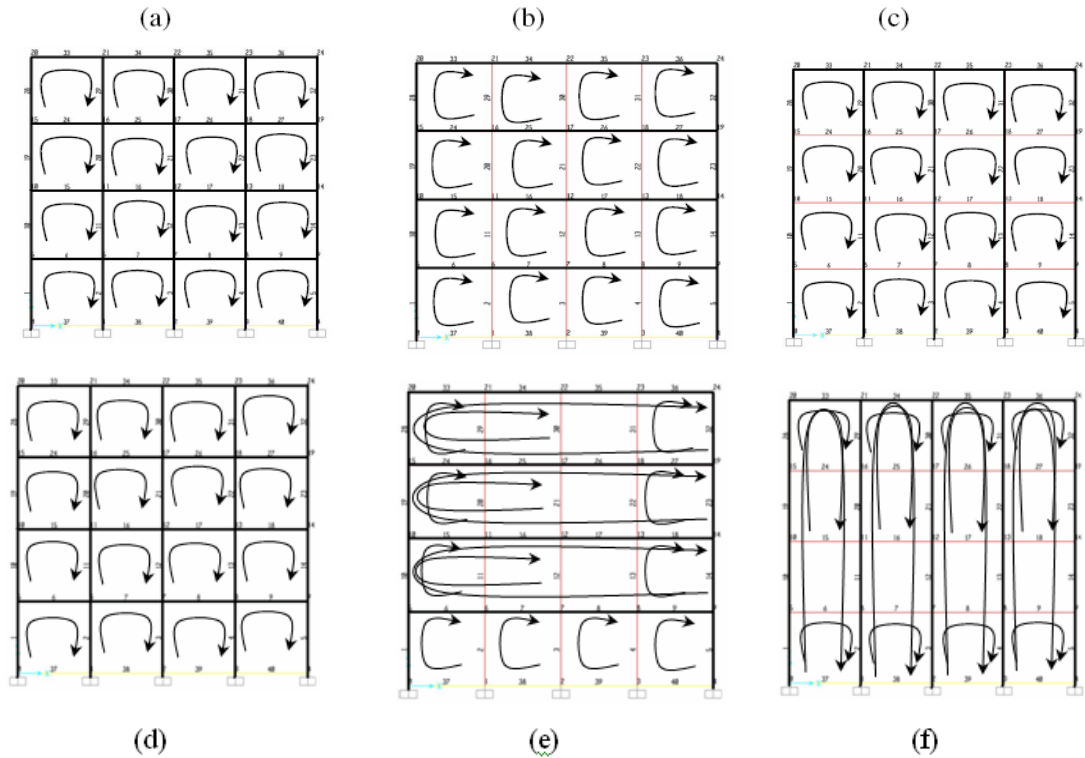

**Figure 11.** 4-story, 4-span frames with different member properties and depicted optimal cycle bases using algorithm 1 (a-c) and algorithm 2 (d-f)

As shown in Figure 11, the frames' optimal cycle bases were calculated using algorithms 1 and 2. The necessary corresponding flexibility matrices for frames are shown in Figure 12(a-c) for algorithm 1 and Figure 12(d-f) for algorithm 2. Due to the fact that both frames (a and d) used similar characteristics for beams and columns, both algorithms produced the same length of cycle bases. Additionally, each cycle is as small as possible, resulting in a more sparse flexibility matrix. Figure 12(a-c).

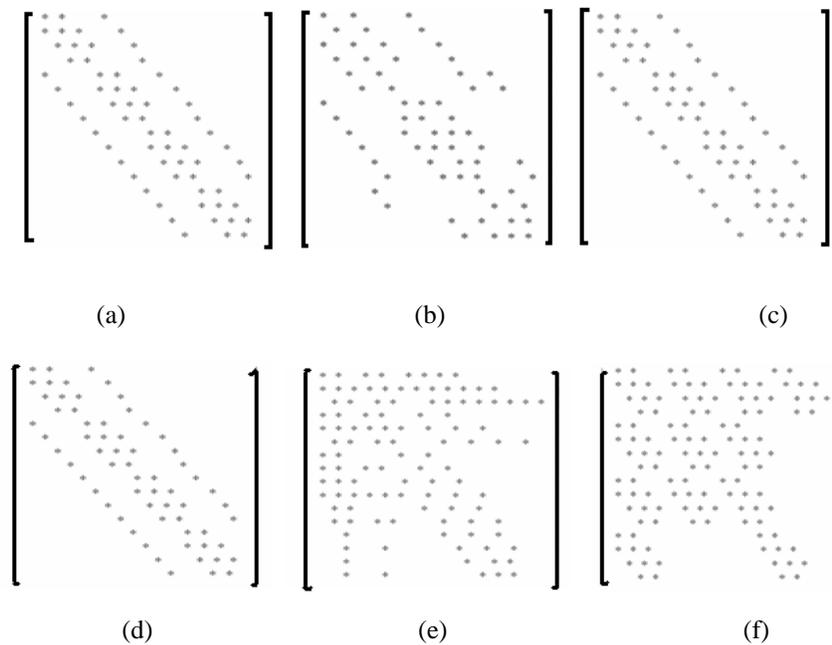

**Figure 12.** Corresponding flexibility matrices for Figure 11 with algorithm 1 (a-c) and algorithm 2 (d-f).

On the contrary, since the properties of beams and columns in the real world are not similar, the beams and columns in the frames were given distinct structural features Figure 11(b-c and e-f). As shown in Figure 12(e-f), the corresponding flexibility matrix is therefore less sparse. When applied to a frame with unique properties, method 1 produces a sparser matrix than algorithm 2. On the other hand, approach 2 generates a less sparse matrix with better well conditioning. When the conditioning number in Table 8 is compared to the conditioning number in Table 9, this truth becomes clearly evident. Additionally, for the frame (d) from Figure 11 , algorithm 2 also generated cycle base with minimum length led to sparser corresponding flexibility matrix due to same assigned structural properties for all beams and columns.

**Table 8.** Well conditioning numbers using algorithm 1

| TYPE | PL | PN | PDET | X(D) |
|---|---|---|---|---|
| (a) | 3.541287 | 0 | 2.296645E-25 | 64 |
| (b) | 4.625147 | 0 | 5.215225E-32 | 64 |
| (c) | 4.479996 | 0 | 0 | 64 |

**Table 9.** Well conditioning numbers using algorithm 2

| TYPE | PL | PN | PDET | X(D) |
|---|---|---|---|---|
| (d) | 3.541287 | 0 | 2.296645E-25 | 64 |
| (e) | 3.958741 | 0 | 1.254628E-29 | 110 |
| (f) | 4.038741 | 0 | 3.837134E-31 | 120 |

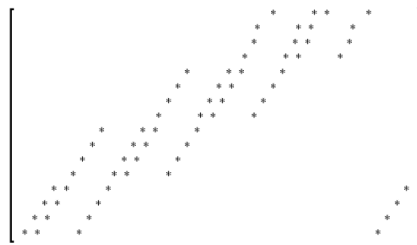

(a)

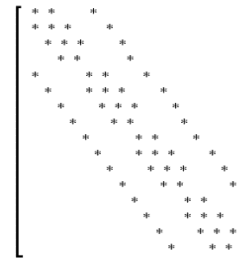

(b)

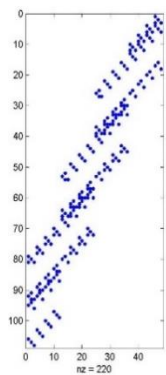

(c)

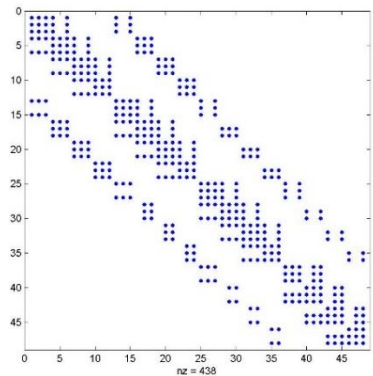

(d)

**Figure 13.** (a) Corresponding self-equation system matrix, (b) Corresponding matrix of flexibility matrix, (c) Matrix $B_1$, and (d) flexibility matrix for model (a) from Figure 11 using algorithm 1

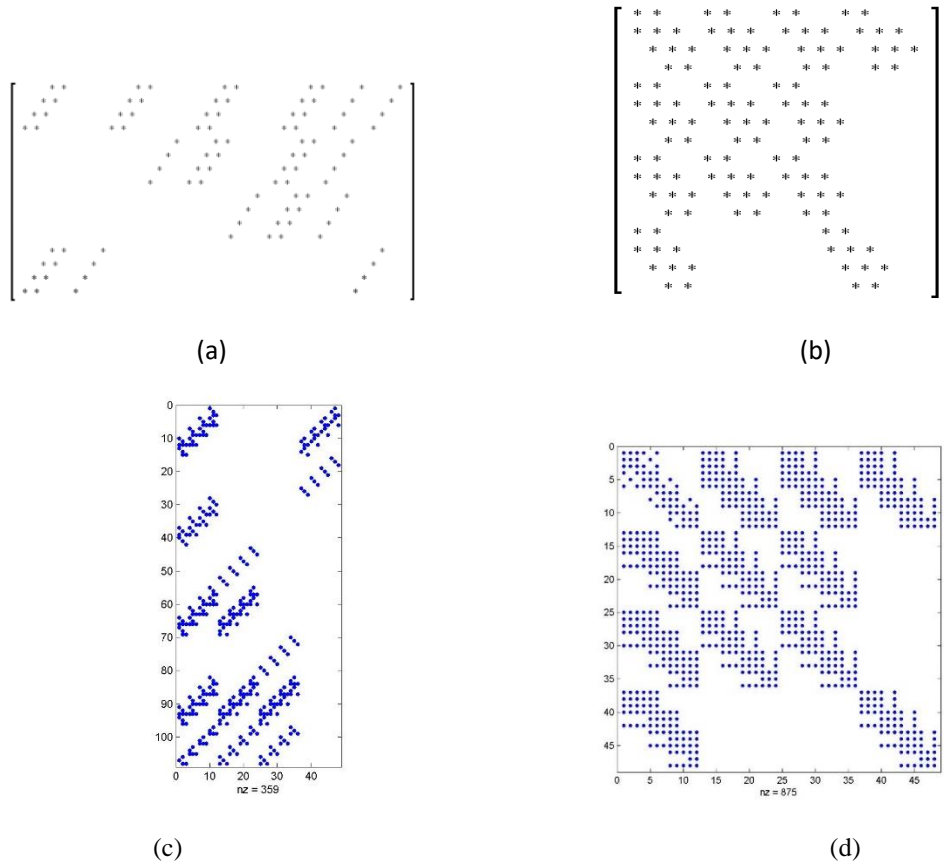

(a)  (b)

(c)  (d)

**Figure 14.** (a) Corresponding self-equation system matrix, (b) Corresponding matrix of flexibility, (c) Matrix $B_1$, and (d) flexibility matrix for model (e) from Figure 11 using algorithm 2

Moreover, Corresponding self-equation system matrices, Corresponding flexibility matrices, Matrices $B_1$, and flexibility matrices for frame model (a) from Figure 11 were shown in Figure 13 using algorithm 1 and also using algorithm 2, same matrices were presented in Figure 14.

**Example 4** – A 9-story frame with members with different properties was shown in Figure 10. In order find sub-optimal fundamental cycle bases algorithms 3 and 4 were applied respectively. Also, the defined properties for beams and columns of the frame were shown in Table 10.

**Table 10.** Defined properties for beams and columns of the frame

| Structural Properties | $A_1$ (m²) | $A_2$ (m²) | $I_1$ (m⁴) | $I_2$ (m⁴) | E (t/m²) |
|---|---|---|---|---|---|
| | 0.00106 | 0.00970 | 0.00000171 | 0.00019610 | $2.1 \times 10^7$ |

The length of the all members were assigned to be 3m.

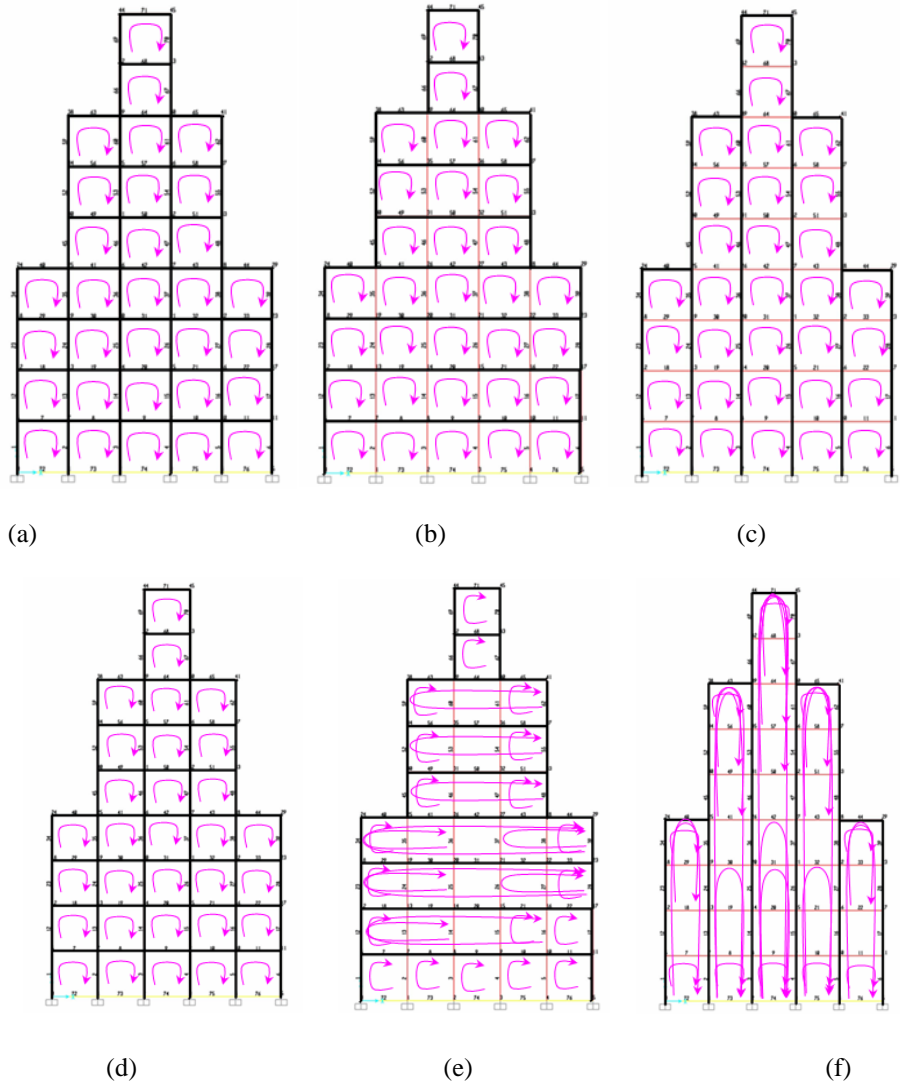

**Figure 15.** 9-story frames with different member properties and depicted optimal cycle bases using algorithm 3 (a-c) and algorithm 4 (d-f)

In Figure 15, two different types of properties were used for beams and columns. As can be seen from the Figure 15, type one of the properties was depicted in the light color, and type two of the properties as shown in the dark color. Additionally, in order to make the supports rigid, members with yellow color were added to the model. Moreover, fundamental cycle bases were obtained using two novel algorithms, and the results were shown and discussed in the following. As shown in Figure 15, the optimum cycle bases for the frames were determined by methods 1 and 2. Figure 15(a-c) depicts the required corresponding flexibility matrices for frames for algorithm 3, and Figure 15(d-f) depicts the required corresponding flexibility matrices for algorithm 4. In light of the fact that both frames (a and d) utilized beams and columns with similar properties, both algorithms generated cycle bases that were the same length. Furthermore, each cycle is kept as short as possible, resulting in a more sparse flexibility matrix. (d)       (e)       (f)

**Figure 16.** Corresponding flexibility matrices for Figure 15 with algorithm 3 (a-c) and algorithm 4 (d-f).

(a-c)

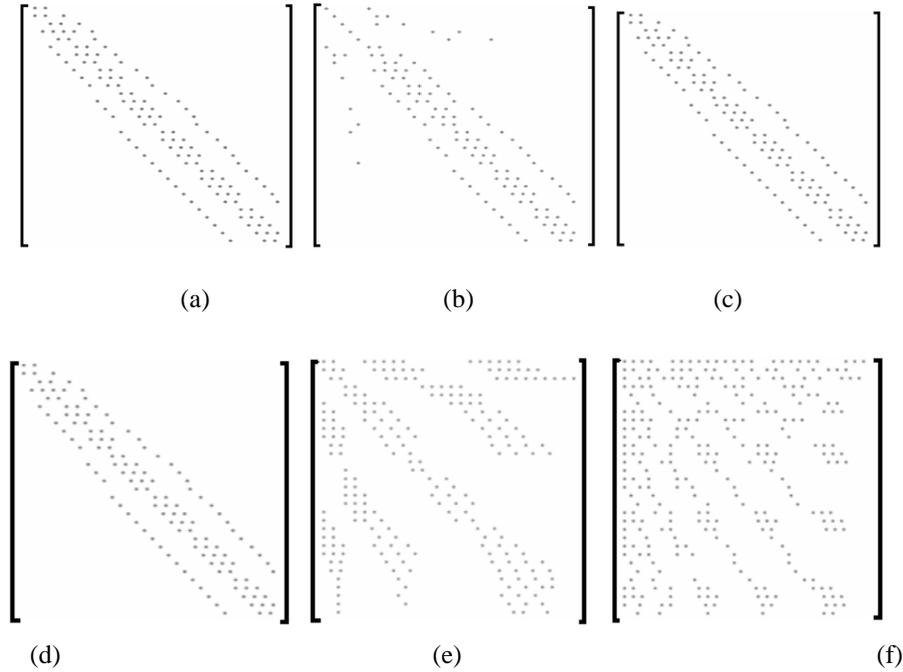

**Figure 16.** Corresponding flexibility matrices for Figure 15 with algorithm 3 (a-c) and algorithm 4 (d-f).

On the contrary, since the features of beams and columns are not identical in the actual world, the beams and columns in the frames were given different structural characteristics Figure 15(b-c and e-f). As shown in Figure 15(e-f), this results in a less sparse flexibility matrix. When applied to a unique frame, the third algorithm generates a sparser matrix than the second algorithm. On the other hand, algorithm 4 produces a matrix that is less sparse and has improved well conditioning. When comparing the conditioning number in Table 11 to the conditioning number in Table 12, this reality becomes abundantly apparent. Additionally, for the frame (d) in Figure 15, algorithm 4 produced a cycle base with the shortest possible length, which resulted in a sparser matching flexibility matrix owing to all beams and columns having the same assigned structural characteristics.

**Table 11.** Well conditioning numbers using algorithm 3

| Type | PL | PN | PDET | X(D) |
|------|----|----|------|------|
| (a) | 4.258306 | 0 | 5.824763E-49 | 127 |
| (b) | 5.125846 | 0 | 1.254628E-62 | 227 |
| (c) | 4.926146 | 3.221330e-73 | 1.486451e-56 | 303 |

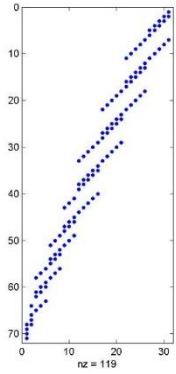
(a)

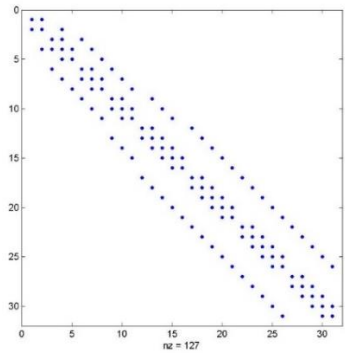
(b)

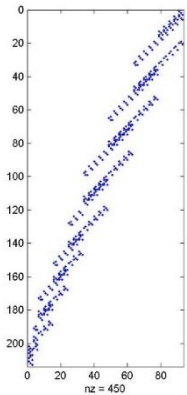
(c)

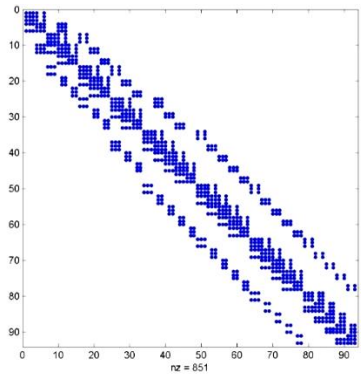
(d)

**Figure 17.** (a) Corresponding self-equation system matrix, (b) Corresponding matrix of flexibility matrix, (c) Matrix B1, and (d) flexibility matrix for model (a) from Figure 15 using algorithm 3

Moreover, Corresponding self-equation system matrices, Corresponding flexibility matrices, Matrices $B_1$, and flexibility matrices for frame model (a) from Figure 15 were shown in Figure 17 using algorithm one and also using algorithm 4; the same matrices were presented in Figure 18. As can be seen from Figure 17 and Figure 18, there is an excellent corresponding between the corresponding flexibility matrix and flexibility matrix itself, which means that the process of matrix optimization could be performed for the corresponding matrix to reduce numerical cost time.

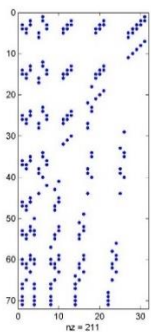
(a)

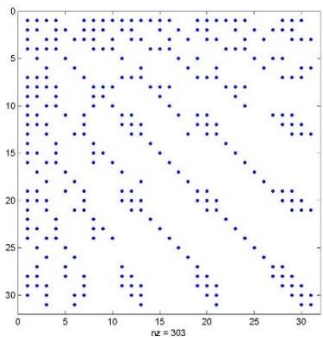
(b)

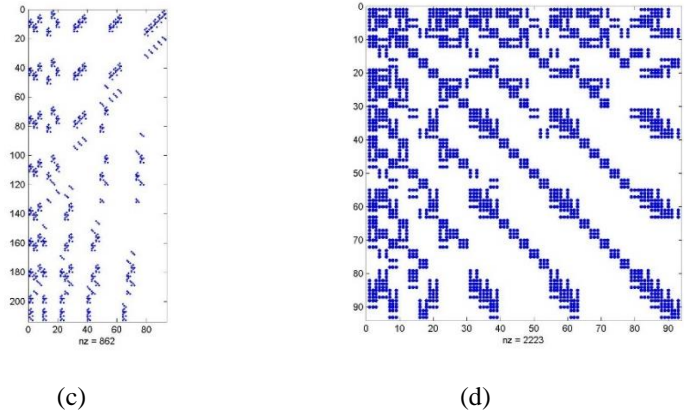

(c)  (d)

**Figure 18.** (a) Corresponding self-equation system matrix, (b) Corresponding matrix of flexibility, (c) Matrix B1, and (d) Flexibility matrix for model (e) from Figure 15 using algorithm 4

**Table 12.** Well conditioning numbers using algorithm 4

| TYPE | PL | PN | PDET | X(D) |
|---|---|---|---|---|
| (a) | 4.258306 | 0 | 5.824763E-49 | 127 |
| (b) | 7.321458 | 0 | 0 | 127 |
| (c) | 6.351422 | 0 | 0 | 127 |

**Example 5** – 4-story space frame with different properties for beams and columns was shown in Figure 19. In order find optimal fundamental cycle bases algorithms 1 and 5 were applied respectively. Two different properties were used for beams and columns shown in Table 13.

**Table 13.** Defined properties for beams and columns of the space frame

| Structural Properties | $A_1$ (m$^2$) | $A_2$ (m$^2$) | $I_1$ (m$^4$) | $I_2$ (m$^4$) | E (t/m$^2$) |
|---|---|---|---|---|---|
| | 0.00106 | 0.00970 | 0.00000171 | 0.00019610 | $2.1 \times 10^7$ |

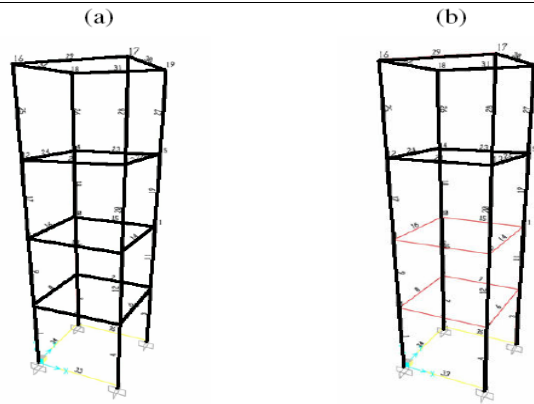

**Figure 19.** A 4-story 3D steel frame with different beams and columns properties.

For the space frame in Figure 19, algorithms 1 and 5 were used to find corresponding flexibility matrix *G*. The results are shown in Figure 12. For frames (a-b), the first algorithm's outcome was presented in Figure 20(a-b). It can be seen that the form of the corresponding flexibility matrix is different, but the sparsity of the two matrices is precisely identical. Additionally, the conditioning number from Table 14 indicates that, well conditioning of matrix b, space frame *b* from Figure 19 with various structural properties for beams and columns, is improved compare to space frame *a* from Figure 19, which dictates that using a frame with the same properties for beams and columns is never recommended, at least for the sake of numerical analysis. The results of the analysis in finding a corresponding flexibility matrix for frames (a-b) using the second algorithm are also shown in Figure 20 (c-d). Comparing sparsity and conditioning of two matrices using Table 14 specifies that frame *a* is sparser than frame *b* but with less improvement in conditioning.

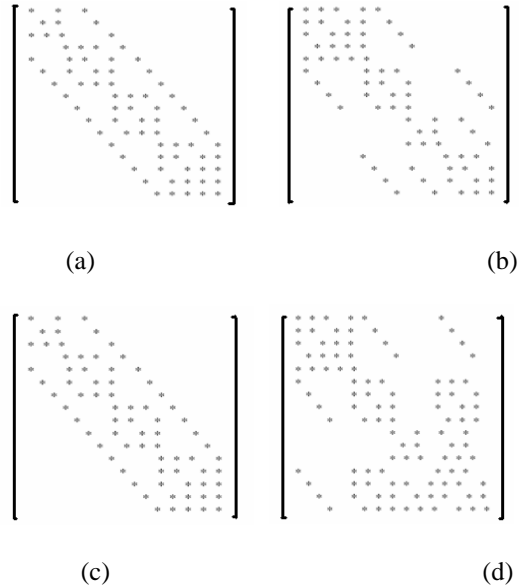

(a)             (b)

(c)             (d)

**Figure 20.** The corresponding flexibility matrices for Figure 18 with algorithm 1 (a-b) and algorithm 5 (c-d).

**Table 14.** Well conditioning numbers using algorithm 1 and 5

| Type | PL | PN | PDET | X(D) |
|---|---|---|---|---|
| Algorithm 1 for model (a) | 6.294568 | 0 | 0 | 74 |
| Algorithm 1 for model (b) | 7.253624 | 0 | 0 | 74 |
| Algorithm 5 for model (c) | 6.294568 | 0 | 0 | 74 |
| Algorithm 5 for model (d) | 5.143547 | 2.213548E-35 | 3284569E-26 | 106 |

## 11. CONCLUSION

In the numerical simulation different kinds of errors may be created during analysis. Round off error is one of the sources of errors. It was shown in the paper that with using different digit decimals for rounding off the final solution, there will be a huge difference in the final answer. Also, using more decimals and rounding off with more decimals is not always a decent method. Because we still do not know how many digit decimals we can select for a specific problem and even if more decimals work for some problems, the computational cost will be high. Dealing with numerical errors in the finite element method is essential. In numerical analysis sometimes handling numerical errors are challenging. However, by applying appropriate algorithms these errors are manageable. In this study novel topological algorithms were proposed in setting up structural flexibility matrix and five different examples were used in applying proposed algorithms.

For any skeletal structure, a graph model can be associated. Then, the problem of optimal structural analysis changes to a problem of combinational optimization. Additionally, having the same pattern between structural matrices and graph model matrices is the essence of switching optimal structural analysis to topological optimization. It was shown in the examples that by applying decent algorithms and associating the relative stiffness of the structural members on the edge of the graph model, the well-conditioning of the structural matrices improved remarkably.

As seen in the examples using well conditioning algorithms decreases the sparsity of the structural matrices, which is not in our favor, but as we discussed, this is a multi-objective problem which means that all three properties of structural matrices cannot be optimized simultaneously. That being said, there should be a compromise between sparsity and the well condition of structural matrices.

## 12. REFERENCES

1. Shah, J.M., Ill-conditioned stiffness matrices. Journal of the Structural Division, 1966. 92(6): p. 443-457.
2. Argyris, J.H., et al., On numerical error in the finite element method. Computer Methods in Applied Mechanics and Engineering, 1976. 7(2): p. 261-282.
3. Grooms, H.R. and J. Rowe, Substructuring and conditioning. Journal of the Structural Division, 1977. 103(3): p. 507-514.
4. Robinson, J. and G.W. Haggenmacher, Some new developments in matrix force analysis. Recent Advances in Matrix Methods of Structural Analysis and Design, Univer. Alabama, 1971: p. 183-228.
5. Henderson, J.d.C., Topological aspects of structural linear analysis: improving the conditioning of the equations of compatibility of a multi-member skeletal structure by use of the knowledge of topology. Aircraft Engineering and Aerospace Technology, 1960.
6. Kaveh, A., Optimizing the conditioning of structural flexibility matrices. Computers & structures, 1991. 41(3): p. 489-494.
7. Maunder, E.A.W., Topological and linear analysis of skeletal structures. 1971.
8. Kaveh, A., The application of topology and metroid theory to the analysis of structures. 1974.
9. Kaveh, A., Improved cycle bases for the flexibility analysis of structures. Computer Methods in Applied Mechanics and Engineering, 1976. 9(3): p. 267-272.
10. Kaveh, A., An efficient program for generating subminimal cycle bases for the flexibility analysis of structures. Communications in applied numerical methods, 1986. 2(4): p. 339-344.
11. Kaveh, A., A combinatorial optimization problem; Optimal generalized cycle bases. Computer Methods in Applied Mechanics and Engineering, 1979. 20(1): p. 39-51.
12. Denke, P.H., A general digital computer analysis of statically indeterminate structures. 1962: National Aeronautics and Space Administration.
13. Robinson, J., Integrated theory of finite element methods(Book on applications to structural and stress analysis). Research supported by the Lockheed-California Co. London and New York, Wiley-Interscience, 1973. 445 p, 1973.
14. Topcu, A., A contribution to the systematic analysis of finite element structures using the force method. Doctoral dissert., Essen Univer, 1979.
15. Kaneko, I., M. Lawo, and G. Thierauf, On computational procedures for the force method. International Journal for Numerical Methods in Engineering, 1982. 18(10): p. 1469-1495.
16. Soyer, E. and A. Topçu, Sparse self-stress matrices for the finite element force method. International Journal for Numerical Methods in Engineering, 2001. 50(9): p. 2175-2194.
17. Gilbert, J.R. and M.T. Heath, Computing a sparse basis for the null space. SIAM Journal on Algebraic Discrete Methods, 1987. 8(3): p. 446-459.
18. Coleman, T.F. and A. Pothen, The null space problem I. Complexity. SIAM Journal on Algebraic Discrete Methods, 1986. 7(4): p. 527-537.
19. Coleman, T.F. and A. Pothen, The null space problem II. Algorithms. SIAM Journal on Algebraic Discrete Methods, 1987. 8(4): p. 544-563.
20. Pothen, A., Sparse null basis computations in structural optimization. Numerische Mathematik, 1989. 55(5): p. 501-519.
21. Patnaik, S., The integrated force method versus the standard force method. Computers & structures, 1986. 22(2): p. 151-163.
22. Patnaik, S., The variational energy formulation for the integrated force method. AIAA journal, 1986. 24(1): p. 129-137.
23. Kaveh, A. and S. Bijari, Simultaneous analysis, design and optimization of trusses via force method. Structural engineering and mechanics: An international journal, 2018. 65(3): p. 233-241.
24. Kaveh, A. and M. Massoudi, Efficient finite element analysis using graph-theoretical force method tetrahedron elements. International Journal of Civil Engineering, 2014. 12(2): p. 249-269.
25. Kaveh, A. and K. Koohestani, Efficient finite element analysis by graph-theoretical force method; triangular and rectangular plate bending elements. Finite elements in analysis and design, 2008. 44(9-10): p. 646-654.
26. Argyris, J.H. and S. Kelsey, Energy theorems and structural analysis. Vol. 60. 1960: Springer.
27. Stepanets, G., Basis systems of vector cycles with extremal properties in graphs. Uspekhi Matematicheskikh Nauk, 1964. 19(2): p. 171-175.
28. Zykov, A., Theory of Finite Graphs [in Russian], Vol. 1, Izd. Nauk (Sibirsk. Otd.), Novosibirsk, 1969.
29. Cassell, A., J. de C. Henderson, and A. Kaveh, Cycle bases for the flexibility analysis of structures. International Journal for Numerical Methods in Engineering, 1974. 8(3): p. 521-528.
30. Hubicka, E. Minimal bases of cycles of a graph. in Recent Advances in Graph Theory. Proc of the symp in Prague (Jun 1974). 1975. Academia Praha.
31. Kolasińska, E., On a minimum cycle basis of a graph. Applicationes Mathematicae, 1980. 4(16): p. 631-639.
32. Kaveh, A. and G. Roosta, Revised Greedy algorithm for formation of a minimal cycle basis of a graph. Communications in numerical methods in engineering, 1994. 10(7): p. 523-530.
33. Horton, J.D., A polynomial-time algorithm to find the shortest cycle basis of a graph. SIAM Journal on Computing, 1987. 16(2): p. 358-366.

34. Kaveh, A. and H. Rahami, Algebraic graph theory for sparse flexibility matrices. Journal of Mathematical Modelling and Algorithms, 2003. 2(2): p. 171-182.
35. Kaveh, A., Structural mechanics: graph and matrix methods. Vol. 6. 1992: Macmillan International Higher Education.
36. DIZAJI, F. and M. Khanzadi, Graph Theoretical Methods for improving The conditioning of Flexibility Matrix of Structures. Journal of Civil Engineering 2011. 22(1): p. 27-44.
37. Dizaji, F.S., Graph Theoretical and Algebra Methods for improving the conditioning of structural Matrices, in Civil Engineering. 2019, Iran University of Science and Technology: Tehran, Iran.